\documentclass[10pt]{article}

\usepackage{theorem,amssymb,amsmath}

\usepackage{graphicx}

\usepackage{color}                    
\definecolor{red}{rgb}{1,0,.2}        
\definecolor{cjp}{rgb}{.1,.7,.2}        
\definecolor{fmdc}{rgb}{1,0,.8}        

\newcommand{\Sh}{{\rm S}}
\newcommand{\Lan}{{\mathcal{ L}}}

\newcommand{\N}{{\mathbb{N}}}
\newcommand{\Card}{{\rm Card}}
\newcommand{\nS}{$\bigtriangledown$}
\newcommand{\aS}{$\boxplus$}
\newcommand{\bS}{$\boxdot$}
\newcommand{\mor}[2]{\Big\{\begin{aligned} a & \rightarrow #1\\[-.1cm] b & \rightarrow #2  \end{aligned}}

\newcommand{\xF}{x_{\rm \scriptstyle F}}
\newcommand{\xG}{x_{\rm \scriptstyle G}}
\newcommand{\xH}{x_{\rm \scriptstyle H}}

\advance \topmargin by -\headheight
\advance \topmargin by -\headsep
\evensidemargin \oddsidemargin
\author{
F. M. Dekking \\
Delft University of Technology \\
Faculty EEMCS, P.O.~Box 5031\\
2600 GA Delft, The Netherlands\\
{\tt F.M.Dekking@math.tudelft.nl}
}

\title{\bf Morphic words, Beatty sequences and integer images of the Fibonacci language}

\date{\today}

\def \proof{\noindent{\it Proof:\ \ }}

\newcommand{\Id}{ {\rm Id}}

\newtheorem{theorem}{Theorem}
\newtheorem{lemma}[theorem]{Lemma}
\newtheorem{corollary}[theorem]{Corollary}
\newtheorem{proposition}[theorem]{Proposition}
\theorembodyfont{\rm}

\newtheorem{remark}[theorem]{Remark}
\newtheorem{example}[theorem]{Example}

\begin{document}

\maketitle

\begin{abstract}
Morphic words are letter-to-letter images of fixed points $x$ of morphisms on finite alphabets. There are situations  where these letter-to-letter maps do not occur naturally, but have to be replaced by a morphism. We call this a decoration of $x$. Theoretically, decorations of morphic words are again morphic words, but in several problems the idea of decorating the fixed point of a morphism is useful. We present two of such problems. The first considers the so called $AA$ sequences, where $\alpha$ is a quadratic irrational, $A$ is the Beatty sequence defined by $A(n)=\lfloor \alpha n\rfloor$, and $AA$ is the sequence $(A(A(n)))$.
The second example considers homomorphic embeddings of the Fibonacci language into the integers, which turns out to lead to generalized Beatty sequences with terms of the form $V(n)=p\lfloor \alpha n\rfloor+qn+r$, where $p,q$ and $r$ are integers.
\end{abstract}

\medskip


\section{Introduction}

A Beatty sequence is a sequence $A = (A(n))_{n \geq 1}$, with $A(n)=\lfloor{n\alpha}\rfloor$ for
$n \geq 1$, where $\alpha$ is a positive real number, and $\lfloor \cdot\rfloor$ denotes the floor function. What Beatty observed is that when
$B = (B(n))_{n \geq 1}$ is the sequence defined by $B(n)=\lfloor{n\beta}\rfloor$, with $\alpha$
and $\beta$ satisfying
\begin{equation}\label{eq:Beat}
\frac1{\alpha}+\frac1{\beta}=1,
\end{equation}
then $A$ and $B$ are \emph{complementary} sequences, that is, the sets
$\{A(n): n \geq 1\}$ and $\{B(n): n\geq 1\}$ are disjoint and their union is the set of positive
integers. In particular if $\alpha = \varphi = \frac{1 + \sqrt{5}}{2}$ is the golden ratio, this gives that
the sequences $(\lfloor n \varphi \rfloor)_{n \geq 1}$ and $(\lfloor n \varphi^2 \rfloor)_{n \geq 1}$
are complementary.

\medskip

It is well known that any composition $U$ of the two sequences $A=(\lfloor n \varphi \rfloor)_{n \geq 1}$ and $B=(\lfloor n \varphi^2 \rfloor)_{n \geq 1}$
 can be written as an integer linear combination $pA+q\Id+r$, where $\Id$ is defined by $\Id(n)=n$, see {\cite[Theorem~13, p.\ 20]{car-sco-hog}}.
 One has for example,
 $$AA=B-1=A+\Id-1,\quad AB=A+B=2A+\Id,\quad BA=A+B-1=2A+\Id-1, \quad BB=A+2B=3A+2\Id.$$

 Such a result does not hold for all quadratic irrationals. I we take, for example, $\alpha=\sqrt{2}$, i.e., we consider the Beatty sequence given by $A(n)=\lfloor n \sqrt{2} \rfloor$, then the complementary Beatty sequence $B$  is given by
 $B(n)=\lfloor n (2+\sqrt{2}) \rfloor.$
 It is proved in \cite{car-Pell} (see also \cite{fraenkel-kim}) that for $n\geq 1$
 $$AB(n)=\lfloor \sqrt{2} \lfloor n (2+\sqrt{2})  \rfloor\rfloor= A(n)+B(n)=2A(n)+2n.$$
However, no expression for $AA$ is given.\footnote{Neither for $BA$. The sequence $BA$ has about the same complexity as $AA$, since $BA=AA+2A$, as implied by $B(n)=A(n)+2n$ for all $n\ge 1$.} In fact, one can prove that there do not exist integers $p,q$ and $r$ such that $AA=pA+q\Id+r$. This follows from Lemma 7 in \cite{GBS} , since the first order difference sequence of $AA$ takes more than 2 values. Still, expressions for $AA$ are known involving  the sequence $\lfloor \sqrt{2}\{n\sqrt{2}\}\rfloor$, see Theorem 1 in \cite{fraenkel1994}, and see \cite{ballot}.
Our Theorem \ref{th:AAGBS} in Section \ref{sec:AA} clarifies the situation.

 We next show in Section \ref{sec:AA} that for an infinite  collection of $\alpha$'s the difference sequence of $AA$, as a word,  can be represented as a decoration of the fixed point of a morphism. We determine this for the Fraenkel family, also known as the \emph{metallic means}. These are the solutions to  $x^2+(t-2)x=t$, where the natural number $t$ is the parameter. For $t=1$ one obtains the golden mean, for $t=2$ the silver mean $\sqrt{2}$.

 Recall that the class of decorations of fixed points of morphisms is equal to the class of morphic words, see, e.g., Corollary 7.7.5 in \cite{all-shall}. In Corollary \ref{cor:morphic} we give the difference sequence of $AA$ for $\alpha=(\sqrt{13}-1)/2$ as a morphic word.

\medskip

In Section \ref{sec:FibL} we present our second example. We solve the Frobenius problem for homomorphic embeddings of the Fibonacci language, which means that we give a precise description of the complement of this embedding. Although the two examples are seemingly unrelated, generalizations of Beatty sequences do appear again.

\medskip

\medskip In the appendix we give a different proof that the iterated Beatty sequence $AA$ with $A(n)=\lfloor n\sqrt{2}\rfloor$ is a morphic word. This leads to a morphic word on an alphabet of size 4. We conjecture that this is the smallest size possible, which is equivalent to the conjecture that $AA$ is not a fixed point of a morphism.

\medskip

For some general results for a special class of decorations of fixed points of morphisms see \cite{frid}. In \cite{frid} the decorations are so called marked morphisms, which in some sense are the opposite of the decorations that one will encounter in the present paper.
We mention also that decorations of morphisms are closely connected to HD0L-systems. See \cite{Seebold-2018} for some recent results on these in the context of Beatty sequences, which in some sense are also opposite to our results.

\section{Iterated Beatty sequences}\label{sec:AA}

Let $\alpha$ be an irrational number larger than 1, and let $A$ defined by $A(n)=\lfloor{n\alpha}\rfloor$ for $n\ge 1$
be the Beatty sequence of $\alpha$.

The iterated Beatty sequence $AA$ given by $AA(n)=\lfloor \lfloor n\alpha\rfloor\alpha\rfloor$ has been studied by many authors. See, among others, \cite{car-sco-hog}, \cite{car-Pell}, \cite{fraenkel1994},  \cite{A-Fraenkel-S}, \cite{ballot}. The main effort in these papers has been to express $AA$ as a linear combination of $A$, $\Id$ and the constant function. Following \cite{GBS} we call any sequence $V$ of the form
$$V(n)=pA(n)+qn+r {\quad\rm for\:} n\ge 1$$
where $p,q,r$ are integers, a {\it generalized Beatty sequence}, for short a GBS.

\medskip

Let $(x-\alpha)(x-\overline{\alpha})$ be the minimal polynomial of a quadratic irrational $\alpha$.

\begin{theorem}  \label{th:AAGBS} Let $\alpha>1$  be a quadratic irrational with minimal polynomial in $\mathbb{Z}[x]$. The sequence $AA$ is a generalized Beatty sequence if and only if \, $|\overline{\alpha}|<1$.
\end{theorem}

\proof
 If one substitutes $K=1$, $L=M=0$, $n=2$ and $a_2=1$ in Theorem 1 of  Fraenkel's 1994 paper \cite{fraenkel1994}, one obtains  $$ AA(n)=-a_1A(n)-a_0n+D(n),               $$
where $(x-\alpha)(x-\overline{\alpha})=x^2+a_1x+a_0$, and (with $\{\cdot\}$ denoting the frational part of a real number)
$$D(n)= \Big\lfloor \frac{a_0}{\alpha}\{n\alpha\} \Big\rfloor.$$
 The theorem now follows, since $\alpha\overline{\alpha}=a_0$, and since the sequence $(\{n\alpha\})$ is  equidistributed over $[0,1]$.\hfill $\Box$

 \bigskip

 \noindent If $S$ is a sequence, we denote its  sequence of first order differences as $\Delta S$, i.e., $\Delta S$ is defined by
$$\Delta S(n) = S(n+1)-S(n), \quad {\rm for\;} n=1,2\dots$$

\noindent {\bf Example} Let $\alpha=1+\sqrt{2}$, with corresponding $A(n)=\lfloor n (1+\sqrt{2}) \rfloor$.
 As in the proof of Theorem \ref{th:AAGBS} one computes  that  $A(A(n))=2A(n) +n-1$. An application of Lemma 7 from \cite{GBS} then gives that the difference sequence $\Delta AA$ given by $\Delta AA(n)= AA(n+1)-AA(n)$ is pure morphic: it is fixed point of the morphism $5\rightarrow 57, \; 7\rightarrow 575$ on the alphabet $\{5,7\}$.

 \bigskip

 \noindent What is the structure of $AA$ if \,$|\overline{\alpha}|>1$? We determine this for the Fraenkel family, also known as the \emph{metallic means}, which are the positive solutions to  $x^2+(t-2)x=t$, where the natural number $t$ is the parameter.

\begin{theorem}\label{thm:metal} Let $\alpha=\big(2-t+\sqrt{t^2+4}\,\big)/2$, for $t=2,3,\dots$, and let $A(n)=\lfloor n \alpha \rfloor$ for $n\ge 1$ . Then $\Delta AA $ is a morphic word. In fact, $\Delta AA $ is a decoration $\delta$ of a fixed point  of a morphism $\tau$, both defined on the alphabet $\{1,2\dots,t\!+1\}$. For $t=2$ and $t=3$ the morphisms $\tau$ and $\delta$ are given respectively by\\[-.3cm]
$$\tau(1)=12, \; \tau(2)= 131,\; \tau(3)=121,\;\; \delta(1)=13,\; \delta(2)=222,\; \delta(3)=132,$$
$$\tau(1)=123, \,  \tau(2)= 124,\, \tau(3)=1141, \, \tau(4)=1241,
\;\;\delta(1)=113,\,  \delta(2)=122, \,  \delta(3)=2122, \, \delta(4)=1222.$$
For $t\ge 4$ the morphism $\tau$ is given\footnote{For readability, we denote the letters $t-j$ as $[t-j]$.} by
$\tau(1) = 1...[t-1]\,t, \;  \tau(2) =1...[t-1]\,[t+1]$,\\ and   for   $j=3,...,t-1$ \\[-.7cm]
\begin{eqnarray*}
  \tau(j) \!\!\! &=& \!\!\! 1...[t-j]\,[t-j+1]\,[t-j+1]\,[t-j+2]\dots [t-2]\,[t+1],\\[-.0cm]
\tau(t) \!\!\! &=& \!\!\!   112\dots [t-2]\,[t+1]\,1,\quad    \tau(t+1) = 1223\dots [t-2]\,[t+1]\,1.
\end{eqnarray*}
For $t\ge 4$ the morphism $\delta$ is given by
$\delta(1) = 1^{t-1}\,3, \;  \delta(2) =1^{t-2}\,22,\;  \delta(j)  =  1^{t-j}\,2\,1^{j-2}\,2$ for $j=3,...,t-\!1,$ \\[.1cm]
and $\delta(t)  =   2\,1^{t-2}\,22, \;    \delta(t+1) = 12\,1^{t-3}\,22.$
\end{theorem}


In the proof of this theorem we need the combinatorial Lemma \ref{lem:comb}.
We know that $\Delta A$ is fixed point of the morphism $\sigma$ on the alphabet $\{1,2\}$ given by
 \begin{equation}\label{eq:sigma} \sigma(1)=1^{t-1}2, \quad \sigma(2)= 1^{t-1}21,\end{equation}
  as can be found in Crisp et al \cite{Crisp-et-al}, or Allouche and Shallit \cite{all-shall}. Here one uses that $\alpha$ has a very simple continued fraction expansion: $\alpha=[1; t,t,t,\dots].$

\begin{lemma}\label{lem:comb}
 Let $t\ge 2$ be an integer. For $t=2$, define the three words  $u_1=121,\;v=2112$, and $w=1212$.

\noindent For $t\ge 3$, define the $t-1$ words  $u_j=1^{t-j}21^i$ for $j=1,...,t-1$, and the two words  $v=21^t2, \; w=121^{t-1}2$.

\noindent Let $\sigma$ be the morphism in {\rm{(\ref{eq:sigma})}}, then for $t=2$, one has $\sigma(u_1) =u_1v, \; \sigma(v) =u_1wu_1, \;\sigma(w) = u_1vu_1$.

\noindent For $t=3$ one has $ \sigma(u_1)=u_1u_2v, \; \sigma(u_2)=u_1u_2w, \; \sigma(v)=u_1u_1wu_1, \; \sigma(w)=u_1u_2wu_1.$

\noindent  For $t\ge 4$ one has
 $\sigma(u_1) = u_1...u_{t-1}v, \;  \sigma(u_2) =u_1...u_{t-1}w$, and   for   $j=3,...,t-1$ one has\\[-.7cm]
  \begin{eqnarray*}
  \sigma(u_j) \!\!\! &=& \!\!\! u_1...u_{t-j} u_{t-j+1}u_{t-j+1}u_{t-j+2}\dots u_{t-2}w,\\[-.0cm]
\sigma(v) \!\!\! &=& \!\!\!   u_1u_1u_2\dots u_{t-2}wu_1,\quad    \sigma(w) = u_1u_2u_2u_3\dots u_{t-2}wu_1.
 \end{eqnarray*}
\end{lemma}

\proof  First we take $t=2$. Then $\sigma$ is given by $\sigma(1)=12,\;  \sigma(2)=121$. One easily verifies the statement of the lemma:
$\sigma(u_1) =1212112= u_1v, \; \sigma(v) =1211212121= u_1wu_1, \;\sigma(w) = 1212112121=u_1vu_1.$\\
The case $t=3$ follows from an analogous computation.

\noindent Next, the case $t\ge 4$. We first mention four relations, directly implied by the definitions,  which will be used in the proof:
  $$v=21\sigma(1), \quad w=12\sigma(1), \quad w=u_{t-1}2, \quad u_1=\sigma(2).$$
We also use repeatedly\\[-.4cm]
 $$\sigma(1^{j})=u_1\dots u_{j-1}1^{t-j}2\quad {\rm for\;} j=2,\dots t-1,$$
 which can be proved by induction: $\sigma(1^{j+1})=u_1\dots u_{j-1}1^{t-j}2\sigma(1)=u_1\dots u_{j-1}1^{t-j}21^{t-1}2=u_1\dots u_{j}1^{t-j-1}2.$

\noindent We then have\\[-.65cm]
 \begin{eqnarray*}
 \sigma(u_1)\!\!\! &=& \!\!\! \sigma(1^{t-1}21)=u_1\dots u_{t-2}12\sigma(2)\sigma(1)=u_1\dots u_{t-2}121^{t-1}21\sigma(1)= u_1\dots u_{t-1}21\sigma(1)=u_1...u_{t-1}v,\\
 \sigma(u_2)\!\!\! &=& \!\!\!\sigma(1^{t-2}211)=u_1\dots u_{t-3}1121^{t-1}211^{t-2}12\sigma(1)=u_1\dots u_{t-1}12\sigma(1)= u_1...u_{t-1}w.
 \end{eqnarray*}
Now for $u_j$, with $3\le j\le t-1$:  (interpreting  $u_1\dots u_{0}$ as an empty prefix in the case $j=t-1$; so in that case the outcome is $\sigma(u_{t-1})=u_1u_{2}u_{2}u_{3}\dots u_{t-2}w$ (if $t\ge 4$).)
\begin{eqnarray*}
 \sigma(u_j) & = & \sigma(1^{t-j})\,\,\sigma(21^j)\\
             & = & u_1\dots u_{t-j-1}1^j2\,\,1^{t-1}21\,\,\sigma(1^j)\\
             & = & u_1\dots u_{t-j-1}\,u_{t-j}\,1^{j-1}21\,\,1^{t-1}2\,\,\sigma(1^{j-1})\\
             & = & u_1\dots u_{t-j-1}u_{t-j}\,u_{t-j+1}\,1^{j-1}2\sigma(1^{j-1})\\
             & = & u_1\dots u_{t-j-1}u_{t-j}u_{t-j+1}\,1^{j-1}2\,\,1^{t-1}2\,\,\sigma(1^{j-2})\\
             & = & u_1\dots u_{t-j-1}u_{t-j}u_{t-j+1}\,u_{t-j+1}\,\,1^{j-2}2\,\,\sigma(1^{j-2})\\
             & = & u_1\dots u_{t-j-1}u_{t-j}u_{t-j+1}u_{t-j+1}\,\,1^{j-2}2\,\,1^{t-1}2\,\,\sigma(1^{j-3})\\
             & = & u_1\dots u_{t-j-1}u_{t-j}u_{t-j+1}u_{t-j+1}\,\,u_{t-j+2}\,\,1^{j-3}2\,\,\sigma(1^{j-3})\\
             & = & \cdots\\
             & = & u_1\dots u_{t-j-1}u_{t-j}u_{t-j+1}u_{t-j+1}\dots u_{t-2}\,\,12\,\,\sigma(1)\\
             & = & u_1\dots u_{t-j} u_{t-j+1}u_{t-j+1}u_{t-j+2}\dots u_{t-2}\,w.
\end{eqnarray*}
For $v$ and $w$ one  derives:
\begin{eqnarray*}
 \sigma(v)\!\!\! &=& \!\!\! \sigma(2)\,\sigma(1^{t-1})\sigma(12)= u_1\,\,u_1\dots u_{t-2}12\,\sigma(1)\,\sigma(2)= u_1u_1u_2\dots u_{t-2}\,w\,u_1\\
 \sigma(w) \!\!\! &=& \!\!\! \sigma(u_{t-1}2)= u_1u_{2}u_{2}u_{3}\dots u_{t-2}w\,
     \sigma(2)=  u_1u_{2}u_{2}u_{3}\dots u_{t-2}w\,u_1. \hspace*{5cm} \Box
\end{eqnarray*}

\bigskip

  {\it{Proof of Theorem \ref{thm:metal}:}}
\; In view of the complexity of the proof we first give the proof for the case $t=3$, i.e., the case $\alpha=(\sqrt{13}-1)/2$, the bronze mean.

 We then have to show that $\Delta AA $ is a decoration $\delta$ of a fixed point  of a morphism $\tau$, both defined on the alphabet $\{1,2,3,4\}$, where $\tau$ is given by\\[-.3cm]
$$\tau(1)=123, \quad  \tau(2)= 124,\quad \tau(3)=1141, \quad \tau(4)=1241.$$
and the decoration  $\delta$ is given by\\[-.3cm]
$$\delta(1)=113,\quad   \delta(2)=122, \quad  \delta(3)=2122, \quad \delta(4)=1222.$$
 The words from Lemma \ref{lem:comb} are in this case\\[-.3cm]
$$ u_1=1121,  \quad u_2=1211,   \quad v=21112, \quad w=12112, $$
and their images under $\sigma$ are\\[-.3cm]
$$   \sigma(u_1)=u_1u_2v, \quad \sigma(u_2)=u_1u_2w, \quad \sigma(v)=u_1u_1wu_1 \quad \sigma(w)=u_1u_2wu_1.$$
The coding $u_1\mapsto 1, u_2\mapsto 2, v\mapsto 3, w\mapsto 4$ transforms $\sigma$ working on $\{u_1,u_2,v,w\}$ into $\tau$.

\medskip

Let $L$ be the map that assigns to any word its length, so, e.g., $L(u_1)=4, L(v)=5$.

\medskip

\noindent CLAIM: 1) The word $\Delta A$ can be written as $\Delta A=x_1x_2\dots$ where each $x_i$ is an element from $\{u_1,u_2,v,w\}$.\\
 2) The word $r:=L(x_1)L(x_2)\dots$ is fixed point of the morphism $\sigma_{4,5}$ given by $4\rightarrow 445, 5\rightarrow 4454$.

\medskip

Proof of part 1) of the claim: we know that $\Delta A$ is the unique fixed point of the morphism $\sigma=\sigma_{1,2}$ given by $1\rightarrow 112, 2\rightarrow 1121$. Since $1121=u_1$ is prefix of $\Delta A$, also $\sigma^n(u_1)$ is prefix of $\Delta A$ for all $n\ge1$. So with Lemma \ref{lem:comb}  this proves the CLAIM, part 1).  Part 2) of the claim then follows from $L(u_1)=L(u_2)=4,\; L(v)=L(w)=5$, which induces the morphism $\sigma_{4,5}$ for the infinite word $r$ of lengths.

\medskip

How do we obtain $\Delta AA$ from $\Delta A$? Since $A(\N)=AA(\N)\cup AB(\N)$, a disjoint union,  one obtains $AA$ from $A$ by removing the integers $AB(n)$, which, of course, have index $B(n)$ in the sequence $A$. The difference sequence $\Delta B$ of this sequence is the unique fixed point of the morphism $\sigma_{4,5}$, since $\beta=\alpha+3$. It follows then from the CLAIM that the integers $AB(n)$ occur at positions with correspond to the third letter in the word $x_i$. Here it is the \emph{third} letter, because the first term of the sequence $(A(B(n))=5,10,15,22,\dots$ occurs at position 4 in the sequence $(A(n))=1,2,3,5,\dots$. Removal of the $AB(n)$ is then performed by adding the third and the fourth letter in the $x_i$. This operation turns $u_1=1121$ into $\delta(1)=1^{t-1}\,3, $, $u_2=1211$ into $\delta(2)=122$, $v=21112$ into $\delta(3)=2122$, and $w=12112$ into $\delta(4)=1222$.
 The conclusion is that this decoration $\delta$ turns the fixed point of $\tau$ into $\Delta AA$.
 This ends the proof for the case $t=3$.

  \medskip

\noindent For general $t$, the coding $u_1\mapsto 1,..., u_{t-1}\mapsto t-1, v\mapsto t, w\mapsto t+1$ transforms $\sigma$ working on $\{u_1,...,u_{t-1},v,w\}$ into $\tau$. An analogous claim as for the $t=3$ case holds, and now  the  map $L$ satisfies
$$L(u_1)=L(u_2)=...=L(u_{t-1})= t+1,\quad L(v)=L(w)=t+2,$$
 which induces the morphism $\sigma_{t+1,t+2}$ for the infinite word $r$ of lengths. One continues in the same way, using now that $\beta=\alpha+t$.
 This time, the integers $AB(n)$ occur at positions in $A$ with correspond to the $t^{\rm th}$ letter in the words $x_i$ from $\{u_1,...,u_{t-1},v,w\}$.  Here it is the \emph{$t^{\rm th}$} letter, because the first term of the sequence $(A(B(n))$ occurs at position $B(1)=t+1$ in the sequence $(A(n))$. Here $B(1)=\lfloor\beta\rfloor=\lfloor\alpha+t\rfloor = t+1$, since a simple computation shows that $1<\alpha<2$ for all $t$.

 Removal of the $AB(n)$ is then performed by adding the  \emph{$t^{\rm th}$} and the  \emph{$(t+1)^{\rm th}$} letter in the $x_i$.
 This operation turns  $u_1=1^{t-1}21$ into $\delta(1)=1^{t-1}3$, $u_2=1^{t-2}211$   into $\delta(2) =1^{t-2}22$ and $u_j=1^{t-j}21^j$ into $\delta(j)=1^{t-j}21^{j-2}2$, for $j=3,...,t-1$. Moreover, the  two words  $v=21^t2, \; w=121^{t-1}2$ are turned into  $\delta(t)  =   2\,1^{t-2}\,22$, respectively    $\delta(t+1) = 12\,1^{t-3}\,22$.

 The conclusion is that this decoration $\delta$ maps the fixed point of $\tau$ to the first differences $\Delta AA$.   \hfill $\Box$

\begin{corollary}\label{cor:morphic} Here is a way to write  $\Delta AA = 11312221222\dots$ as a morphic word for the case $t=3$.
Let $\theta$ on $\{1,\dots, 6\}$ be the morphism given by\\
\hspace*{1cm} $\theta: \quad 1\rightarrow 123,\; 2\rightarrow 164, \;  3\rightarrow  5145, \; 4\rightarrow 1645, \;5\rightarrow 123, \; 6\rightarrow 164.$\\
Let the letter-to-letter morphism $\lambda$ be given by\\
\hspace*{1cm} $\lambda: \quad  1\rightarrow 1,\; 2\rightarrow 1,\; 4\rightarrow 2,\; 5\rightarrow 2,\; 6\rightarrow 2,\; 3\rightarrow 3.$\\
Then  $\Delta AA = \lambda( \theta^{\infty}(1)).$
\end{corollary}

\noindent Corollary \ref{cor:morphic} is derived from Theorem \ref{th:AAGBS} by  using the natural algorithm given, for example, in \cite{Honkala}, Lemma 4. Honkala's requirement of `cyclicity' in that lemma is not necessary.

\bigskip

 Fraenkel's theorem with the 'defect' function $D=D(n)$ suggests that the $\Delta AA$ sequences can take many values. This is not the case.

 \medskip

 \begin{proposition} For any  irrational $\alpha$ larger than $1$ the sequence
 $\Delta AA=\big( \lfloor \lfloor (n+1)\alpha\rfloor\alpha\rfloor\,-\,\lfloor \lfloor n\alpha\rfloor\alpha\rfloor\big)$ takes values in an alphabet of size two, three or four.
 \end{proposition}

 \medskip

 \proof We illustrate the proof with the case $1<\alpha<2$. Then $s:=\Delta A$ is a Sturmian word taking values $d=1$ or $d=2$. So
 $$\Delta AA(n)= A(A(n+1))-A(A(n))= A(A(n)+d)-A(A(n)), \;{\rm where\:} d=1\;{\rm or\;} 2.$$
We put $i:=A(n)$. In case $d=1$, $A(A(n)+d)-A(A(n))=A(i+1)-A(i)=1$ or 2. In case $d=2$,  $A(A(n)+d)-A(A(n))=A(i+2)-A(i)=A(i+2)-A(i+1)+A(i+1)-A(i)$.
So either  $A(A(n)+d)-A(A(n))=2$ or 3, or $A(A(n)+d)-A(A(n))=3$ or 4,  respectively if 11, 12 and 21 are the subwords of length 2 of $s$, or if 12, 21 and 22 are the subwords of length 2 of $s$. What we found is that $\Delta AA$ takes values in $\{1,2,3\}$ if $1<\alpha<3/2$, and $\Delta AA$ takes values in $\{1,2,3,4\}$ if $3/2<\alpha<2$.
 In some cases $\Delta AA$ may take only 2 values, for example, if $\alpha$  is the golden ratio.\\
The proof for other values of $\alpha$ is similar, exploiting balancedness of the Sturmian word $s=\Delta A$.
  \hfill $\Box$

  \bigskip

\noindent   {\bf  Example} Take $\alpha = \sqrt{11}/2 = 1.658\dots$. Then $ AA(n)=1, 4, 6, 9, 13, 14, 18, 21, 23\dots$, so $\Delta AA$ takes the four values $1,2,3$ and $4$.

\bigskip

\noindent   {\bf  Remark} Once more, let $\alpha=\sqrt{2}$. The differences $x_{2,k}:=(AB)^k-(BA)^k$, where $k\ge 1$, are the `commutator' functions. They are extensively studied in \cite{Chew-Tanny}. They are all similar to $x_{2,1}$, which is equal to $x_{2,1}=AB-BA=2\Id-AA$. One can derive from this that all commutator functions are morphic words.

\section{Embeddings of the Fibonacci language into the integers}\label{sec:FibL}

Let $\Lan$ be a language, i.e., a sub-semigroup of the free semigroup generated by a finite alphabet
under the concatenation operation.
A homomorphism of $\Lan$ into the natural numbers is a map $\Sh:\Lan\rightarrow {\mathbb N}$
satisfying $$\Sh(vw)=\Sh(v)+\Sh(w),\quad{\rm  for\: all\;}  v,w \in \Lan.$$

Let $\xF$ be the Fibonacci word, i.e., the infinite word fixed by the morphism $0\rightarrow 01, \,1\rightarrow 0$. Let $\Lan_{\rm \scriptstyle F}$ be the Fibonacci language, i.e., the set of all words occurring in
$x_{\rm \scriptstyle F}$. Recall that $\varphi=(1+\sqrt{5})/2$. The key ingredient in this section is  the lower Wythoff sequence $(\lfloor n\varphi \rfloor)_{n \ge 1}= 1, 3, 4, 6, 8, 9, 11, 12, 14, 16, 17, 19,\dots$. The following result is proved in \cite{Dekking-TCS-2018}.

 \begin{theorem}\label{th:Fib}{\bf \rm (\cite{Dekking-TCS-2018})} Let  $\Sh: \Lan_{\rm \scriptstyle F}\rightarrow \mathbb{N}$ be a homomorphism. Define $a=\Sh(0), b=\Sh(1)$. Then  $\Sh(\Lan_{\rm \scriptstyle F}$) is the union of the two generalized Beatty sequences   $ \big((a-b)\lfloor n\varphi \rfloor+(2b-a)n\big)$ and $ \big((a-b)\lfloor n\varphi \rfloor+(2b-a)n+a-b\big)$.
\end{theorem}

 The goal of this section is to determine the complement of the set $\Sh(\Lan_{\rm \scriptstyle F}$) in $\mathbb{N}$. We shall show that the corresponding infinite word is always a morphic word, by representing it as a decoration of a fixed point of a morphism. It appears that this is a matter of a complicated bookkeeping, especially when  the two values $\Sh(0)$ and $\Sh(1)$ are small.

 There are three morphisms $f,g$ and $h$ that play an important role in this section, where it is convenient to look at $a$ and $b$ both as integers and as abstract letters. The morphisms are given by
$$f:\mor{ab}{a}\,,\qquad g:\mor{baa}{ba}\,,\qquad h:\mor{aab}{ab}\,.$$

\begin{lemma} Let $\xF$ be the Fibonacci sequence  on the alphabet $\{a,b\}$, fixed point of $f$. Then the fixed point $\xG$ of $g$ is the  sequence  $b\,\xF$, and the fixed point $\xH$ of $h$ is $a\,\xF$.
\end{lemma}

\proof See Theorem 3.1 and Example 1  in \cite{Berstel-Seebold}. \hfill $\Box$

\bigskip

\noindent Here is a result that gives an idea of  the proof in general for the case $\Sh(0)> \Sh(1)$.


\begin{theorem}\label{th:dec1} Let  $\Sh: \Lan_{\rm \scriptstyle F}\rightarrow \mathbb{N}$ be a homomorphism determined by $a=\Sh(0), b=\Sh(1)$. Suppose that
$$a+2<2b+1<2a-1.$$ Then the first differences of the complement $\mathbb{N}\setminus\Sh(\Lan_{\rm \scriptstyle F}$) of $\Sh(\Lan_{\rm \scriptstyle F}$) is the word obtained by decorating the fixed point $\xH$ of the morphism $h$ by the morphism $\delta$ given by
$$\delta(a)=1^{b-2}\,2\,1^{a-b-2}\,2, \quad \delta(b)=1^{2b-a-2}\,2\,1^{a-b-2}\,2.$$
\end{theorem}

\proof The sequence of first differences of a generalized Beatty sequence   $ \big(p\lfloor n\varphi\rfloor+qn+r)$ is the fixed point of the Fibonacci morphism $f$  on the alphabet $\{2p+q, p+q\}$. See Lemma 8 in \cite{GBS}. So the two generalized Beatty sequences
$G_1:=\big((a-b)\lfloor n\varphi \rfloor+(2b-a)n\big)$ and $G_2$, given by $G_2(n)=G_1(n)+a-b$ in Theorem~\ref{th:Fib} have the property that $\Delta G_1=\Delta G_2$ is the fixed point $\xF$ of the Fibonacci morphism on the alphabet with symbols $2(a-b)+2b-a=a$ and $a-b+2b-a=b$.

We illustrate the proof by first considering the case $a=8,\,b=5$.
In this case we have
$$G_1=5,13,18,26,34,39,47,52,60,\dots, \quad G_2=G_1+3=8, 16, 21, 29, 37, 42, 50, 55, 63,\dots...$$

\noindent Partition the positive integers $\N$ into adjacent sets $V_i, \, i=1,2,\dots$ defined by
$$V_i=\{G_2(i-1)+1,\dots,G_2(i)\}.$$
Here we put $G_2(0)=0$.
As a consequence, $\Card(V_i) = 8$ if $\xH(i)=a$ and $\Card(V_i)=5$ if $\xH(i)=b$,
where $\xH=\xH(1)\xH(2)\dots=a\,a\,b\,a\,a\,b\dots$ is the fixed point of $h$.
The reason that the directive sequence is $\xH$ instead of $\xF$ is that the \emph{last} element of each $V_i$ is equal to  $G_2(i)$ for $i=1,2,\dots$.

\bigskip

{\footnotesize
\hspace*{-.9cm} \begin{tabular}{ cccccccc|cccccccc|ccccc|ccc }
 &  &  & $V_1$ &  &  &  &  &  &  &  & $V_2$ &  &  &  &  &  &  & $V_3$ &  &  &  &  &  \\
  1 & 2 & 3 & 4 & 5 & 6 & 7 & 8 & 9 & 10 & 11 & 12 & 13 & 14 & 15 & 16 & 17 & 18 & 19 & 20 & 21 & 22 & 23 & 24\\
  \nS & \nS & \nS & \nS & \bS & \nS & \nS & \aS &\nS &\nS & \nS & \nS & \bS & \nS& \nS& \aS& \nS &\bS  &\nS & \nS&\aS & \nS &\nS & \nS\\
\end{tabular} }

\bigskip

\noindent In the table above, the integers in $G_1(\N)$ are marked with  \bS, those in $G_2(\N)$ with  \aS, and those in the complement with a \nS.
By construction, \emph{all} the $V_i$ with cardinality 8 have the same pattern   \nS\nS\nS\nS\bS\nS\nS\aS\: for their members. Also \emph{all} $V_i$ with cardinality 5 have the same pattern \nS\bS\nS\nS\aS. Note that the last two symbols are \nS\aS, for both size 5 and size 8 $V_i$'s, \emph{and} their first symbols are \nS\, for both. This implies that if we glue the patterns together, then the infinite sequence of differences of the positions of \nS\: in the infinite pattern yields first differences of the sequence of elements in $\N\setminus (G_1(\N)\cup G_2(\N))$. For  $V_i$ of size 8 these differences (including the 'jump over' last value 2) are given by 1,1,1,2,1,2, and  for $V_i$ of size 5 by 2,1,2. It follows that the first differences are obtained by decorating the fixed point $\xH$ by the morphism $\delta$ given by
$$\delta:\quad \;a\rightarrow 111212, \quad b\rightarrow 212.$$
For the general case one considers sets $V_i$ of consecutive integers of size $a$ or size $b$, where the order is again dictated by the fixed point $\xH$ of $h$.
The corresponding patterns have exactly one symbol \aS\, at the end, and exactly one symbol \bS\, positioned $a-b$ places before the end. It follows again that over the $V_i$'s the first differences of the complement set end in 2 (the `jump over' value), are preceded by $a-b-2$ 1's, which is preceded by a 2. The first differences start with a number of 1's, which is $(a-2)-1-(a-b-2)-1=b-2$ for the $V_i$'s of length $a$, and $(b-2)-1-(a-b-2)-1=2b-a-2$ for the $V_i$'s of length $b$. This yields the decoration $\delta$ stated in the theorem. \hfill $\Box$

\bigskip

We now give an example of the difficulties one encounters when $\Sh(0)$ or  $\Sh(1)$ are (relatively) small.

\begin{theorem}\label{th:dec31} Let  $\Sh: \Lan_{\rm \scriptstyle F}\rightarrow \mathbb{N}$ be the homomorphism determined by $a=\Sh(0)=3,\, b=\Sh(1)=1$.  Then the sequence of first differences of the complement $\N\setminus\Sh(\Lan_{\rm \scriptstyle F})$ of $\Sh(\Lan_{\rm \scriptstyle F})$ is the word obtained by decorating the fixed point $\xH$ of  $h$ by $\delta:\{a,b\}\rightarrow\{7,11\}$ given by $\,\delta(a)=7,11$, and $\,\delta(b)=11.$
\end{theorem}

\proof According to Theorem \ref{th:Fib}, $\Sh(\Lan_{\rm \scriptstyle F})$ is the union of the two sets $G_1(\N)$ and $G_2(\N)$ given by
$$G_1(\N)=\{2\lfloor n\varphi \rfloor-n,\, n\ge 1\}=1,4,5,8,11,12\dots, \; G_2(\N)=\{2\lfloor n\varphi \rfloor-n+2,\, n\ge 1\}=3,6,7,10,13,14\dots.$$
The first differences $\Delta G_1=\Delta G_2$ are the Fibonacci word on the alphabet $\{3,1\}$. Imitating the proof of the previous theorem, we obtain  the following table, induced by the morphism $h$ given by $1\rightarrow 331, 3\rightarrow 31$.
 One has $\Card(V_i) = a=3$ if $\xH(i)=a$ and $\Card(V_i)=b=1$ if $\xH(i)=b$, where $\xH=\xH(1)\xH(2)\dots=a\,a\,b\,a\,a\,b\,a\,b\,a\,a\dots$ is the fixed point of $h$.

\bigskip

{\footnotesize
\hspace*{-.9cm} \begin{tabular}{ ccc|ccc|c|ccc|ccc|c|ccc|c|ccc|ccc }
    &$V_1$&   &   &$V_2$ &  &$V_3$&  &$V_4$&   &    &$V_5$&   &$V_6$&   &$V_7$&  &$V_8$&  &$V_9$&  &  &$V_{10}$  &  \\
  1 &   2 & 3 & 4 & 5    & 6& 7   & 8&   9 & 10& 11 & 12  & 13& 14  & 15& 16  &17& 18  &19& 20  &21&22& 23       & 24\\
\bS &\nS&\aS&\bS&\bS&\aS&\aS&\bS&\nS&\aS&\bS &\bS &\aS &\aS& \bS& \bS& \aS& \aS&\bS &\nS&\aS &\bS &\bS &\aS\\
\end{tabular} }

\bigskip

\noindent There are at least two things wrong with this:
\begin{enumerate}
\item[[E1\!\!]] The $V_i$'s of length 3 do not all have the same pattern,
\item[[E2\!\!]] There are patterns that do not contain a \nS.
\end{enumerate}
To counter these problems, we go from the letters $a=3,b=1$ to the words $h(3),h(1)$, yielding a partition with  $W_i$'s of length 7 and 4. The table we obtain is

\bigskip

{\footnotesize
\hspace*{-.9cm} \begin{tabular}{ ccccccc|ccccccc|cccc|cccccc }
&   &   &  $W_1$ &&   &  &  &   &   &  $W_2$ &    &  &&  &$W_3$     &   &   &   &  &  &$W_4$    & &      \\
1& 2 & 3 & 4 & 5   &6  & 7& 8& 9 & 10& 11& 12 &13& 14  &15& 16 &17  & 18 &19  & 20&21 &22&23 & 24 \\
\bS&\nS&\aS&\bS&\bS  &\aS&\aS&\bS&\nS&\aS&\bS &\bS &\aS &\aS& \bS& \bS& \aS& \aS&\bS &\nS&\aS&\bS&\bS&\aS \\
\end{tabular} }

\bigskip

\noindent Problem [E1] is caused by the fact that $V_i$'s of length 3 have different patterns depending on whether they are followed by a $V_i$ of length 1 or of length 3. Problem [E1] is now solved with the $W_i$'s, since 33 can only occur as a prefix of $h(1)=331$, and 31 can only occur as a suffix of either $h(1)$ or $h(3)$.

\noindent However,  [E2] is not yet solved, since $W_3$ does not contain a \nS. The way to tackle this is to pass to the square of $h$, i.e., take the $W'_i$'s of length 18 and 11 corresponding to $h^2(1)=33133131$ and $h^2(3)=33131$.

It is obvious from the corresponding patterns, that the differences of the complement $\N\setminus\Sh(\Lan_{\rm \scriptstyle F}$) are given by the decoration $W_1'\rightarrow 7,11, W_3'\rightarrow 11$ of the $W_i'$'s. But since $h^2(\xH)=\xH$, this is the same as decorating the letters $a\rightarrow 7,11$, and $b\rightarrow 11$ in $\xH$.\hfill $\Box$

\begin{remark}\label{rem:GBS} Theorem \ref{th:dec31} gives essentially the same result as Theorem 25 in \cite{GBS}. The proof given here is completely different.
\end{remark}


We let $C$ be the increasing sequence of integers in the complement of $\Sh(\Lan_{\rm \scriptstyle F})$, so $C(\N)=\N\setminus \Sh(\Lan_{\rm \scriptstyle F})$.

\begin{theorem}\label{th:fibL}  Let  $\Sh: \Lan_{\rm \scriptstyle F}\rightarrow \mathbb{N}$ be a homomorphism.  Then the sequence $\Delta C$ of first differences of the complement \,$\N\setminus\Sh(\Lan_{\rm \scriptstyle F}$) of \,$\Sh(\Lan_{\rm \scriptstyle F}$) is a fixed point of a morphism on an alphabet of two letters decorated by a morphism $\delta.$
\end{theorem}

\proof The homomorphism $\Sh$ is determined by $a:=\Sh(0), b:=\Sh(1)$.

\medskip

Case 1: $a\ge 4, b=1$. Here we follow the proof of Theorem \ref{th:dec31}. The $V_i$ are given by $V_i=\{G_2(i-1)+1,\dots,G_2(i)\}.$ Problem [E1], mentioned in the proof of Theorem \ref{th:dec31}, is more severe in this case, as the pattern of the $V_i$'s of length $a$ depends both on $V_{i-1}$ and $V_{i+1}$. If these have both length 1, then the distance to the next element in $V_{i+1}$ with symbol \nS\: is 5, otherwise it is 4. To make the process context free, we  choose the $W_i$ corresponding to the two words
$$v:=h^2(a)=aa1aa1a1, \quad w:=h^2(1)=aa1a1.$$
Context-freeness now occurs because $1a1$ occurs uniquely inside $v$ and $w$. One checks that the
  decoration $\delta$ is then given by
$$\delta(v)=1^{a-3}\,4\,1^{a-4}\,4\,1^{a-3}\,4\,1^{a-4}\,5\,1^{a-4}\,4,\quad \delta(w)=1^{a-3}\,4\,1^{a-4}\,5\,1^{a-4}\,4.$$
Since $v$ and $w$ start and end with the same words, this decoration yields $\Delta C$, when applied to $\xH$ on the alphabet $\{v,w\}$.

\medskip

Case 2: $a=1, b\ge 5$. This\footnote{We leave the case $a=1, b=4$ as an exercise to the reader. In this case the decoration $\delta$ turns out to be $v\rightarrow 11, w\rightarrow 17$.}  is a variant of Case 1. The sequence $\Delta G_1$ is fixed point of the Fibonacci morphism on the alphabet $\{1,b\}$, and so $b\,\Delta G_1$ is fixed point of $g$ on $\{1,b\}$. Problem [E1] is now that the 'jump over' from $b11$ to $b1b$ is 6, but the `jump over' from $b11$ to $b11$ equals 7. The adequate partition elements $W_i$  correspond to the words  $v$ or $w$:
$$v:=g^2(1)=b\,1\,b\,1\,1\,b\,1\,1, \quad  w:=g^2(b)=b\,1\,b\,1\,1.$$
\noindent The decoration $\delta$ is given by
$$\delta(v)=1^{b-4}\,6\,1^{b-5}\,7\,1^{b-5}\,6,\quad \delta(w)=1^{b-4}\,6\,1^{b-5}6.$$

\medskip

Case 3: $a>b\ge 2 $. The partition elements are defined as $V_i=\{G_2(i-1)+1,\dots,G_2(i)\}$, where we put $G_2(0)=0$. This gives  $\Card(V_i) = a$ if $\xH(i)=a$ and $\Card(V_i)=b$ if $\xH(i)=b$, where $\xH=\xH(1)\xH(2)\dots=a\,a\,b\,a\,a\,b\dots$ is the fixed point of $h$. To get rid of problem [E1], we coarsen the partition to blocks $W_i$  corresponding to the words $h(a)=aab$ and $h(b)=ab$.
The problem disappears because $aa$ uniquely occurs as a prefix of $aab$, and $ab$ uniquely as a suffix of $aab$ or $ab$.  Problem [E2] will not occur, since any 5 consecutive integers will contain an element of $C$ (as $b\ge 2$, and no $bb$ occurs in $\xH$), and the smallest cardinality of a $W_i$ is $a+b\ge 5$.
 Also, since both $aab$ and $ab$ start with $a$, and both end in $b$, the patterns of the $W_i$ will concatenate consistently, so that the the decoration $\delta$ obtained form the patterns of the $W_i$  acting as a morphism on $\xH$, will yield the difference sequence of $C$.

\medskip

Case 4: $b>a\ge 2 $. The partition elements are defined as $V_i=\{G_1(i-1)+1,\dots,G_1(i)\}$, where we put $G_1(0)=0$. This gives  $\Card(V_i) = a$ if $\xG(i)=a$ and $\Card(V_i)=b$ if $\xG(i)=b$, where $\xG=\xG(1)\xG(2)\dots=b\,a\,b\,a\,a\,b\dots$ is the fixed point of $g$.
The rest of the proof follows Case 3, replacing $h$ by $g$ (noting that this time $aa$ uniquely occurs as a \emph{suffix} of $g(a)$, and $ab$ only occurs split over a suffix of $g(\ell)$ and a prefix of $g(\ell')$, for $\ell,\ell' = a,b$).\hfill $\Box$

\bigskip

\noindent We illustrate Case 4 with the following example.

\begin{example}\label{ex:59} Let $a=5$, $b=9$. Then $G_1=9,14,23,\dots$ and $G_2=5,10,19,\dots$. The partition elements are $W_1$ of cardinality 14 corresponding to $g(b)=ba=95$, and $W_2$ of cardinality 19, corresponding to  $g(a)=baa=955$.

The patterns of these sets are
\nS\nS\nS\nS\aS\nS\nS\nS\bS\aS\nS\nS\nS\bS\: and  \nS\nS\nS\nS\aS\nS\nS\nS\bS\aS\nS\nS\nS\bS\aS\nS\nS\nS\bS.\\
It follows that the decoration $\delta$ is given by $\delta(9)=\delta(b)=1112113112, \; \delta(5)=\delta(a)=1112113113112.$\hfill $\Box$
\end{example}

\medskip

\noindent The representation in Theorem \ref{th:fibL} is by no means unique.
As an example, let the morphism $\hat{g}_2$ on $\{1,2,3\}$ be given by\\[-.1cm]
\hspace*{4cm} $\hat{g}_2(1)=12,\; \hat{g}_2(2)=\hat{g}_2(3)=132.$\\[.1cm]
\noindent The morphism $\hat{g}_2$ is the 2-block morphism of $g$  under the coding $ba\rightarrow 1, ab\rightarrow 2, aa\rightarrow 3$
 (cf.~\cite{Queff} and \cite{Dekking-JIS}).
The use of $\hat{g}_2$ gives an alternative way to solve problem [E2], leading, for example, in Example \ref{ex:59} to the fact that $\Delta C$ is the decoration of the fixed point of  $\hat{g}_2$ by the morphism $\delta$ given by\\[.1cm]
\hspace*{4cm} $\delta(1)=1112113,\: \delta(2)=112,\: \delta(3)=113.$

\bigskip

Finally we mention another way in which the representation in Theorem \ref{th:fibL} is not unique. In fact, one can show that every $\Delta C$ is a decoration of the single word $\xG$. Let $\bar{f}$ be the time reversal of the Fibonacci morphism $f$, i.e., $\bar{f}$ is defind by $\bar{f}(0)=10,\:\bar{f}(1)=0$. One verifies that \\[.1cm]
\hspace*{4cm} $g=\bar{f} f, \quad h=f \bar{f}.$\\[.1cm]
This leads to \\[.1cm]
\hspace*{4cm} $\bar{f}(\xH)=\bar{f}h(\xH)=\bar{f}f \bar{f}(\xH)=g\bar{f}(\xH)\quad \Rightarrow \quad \bar{f}(\xH)=\xG,$\\[.1cm]
since $\xG$ is the unique fixed point of $g$.
As a corollary one obtains that if $z$ is a decoration of $\xH$ by $\delta$, then $z$ is \emph{also} a decoration of $\xG$: replace $\delta$ by $\delta'=\bar{f}\delta$.

This relation between $\xH$ and $\xG$ is also useful in establishing the connection mentioned in Remark \ref{rem:GBS}.

\medskip

\section{Appendix}
In this section we give an alternative proof of Theorem \ref{thm:metal}, when $t=2$, i.e., the case $\alpha=\sqrt{2}$.

\medskip

\begin{theorem}\label{thm:sol} Let $\alpha=\sqrt{2}$, $A(n)=\lfloor n \alpha \rfloor$ for $n\ge 1$ . Then $\Delta AA = 1,3,2,2,2,1,3,1,3,2,\dots$ is a decoration $\delta$ of a fixed point  of a morphism $\sigma$, both defined on the alphabet $\{1,2,3\}$. Here $\sigma$ is given by\\[-.3cm]
$$\sigma(1)=123, \quad  \sigma(2)= 1,\quad \sigma(3)=121,$$
\vspace*{-.5cm}and the decoration  $\delta$ is given by\\[-.1cm]
$$\delta(1)=13,\quad   \delta(2)=2, \quad     \delta(3)=22.$$
\end{theorem}

 \proof {  Step 1.} In this step we `refine' the sequence $x:=\Delta A$ to a sequence $y$ on 4 symbols, which codes the occurrence of the terms of $AA$ in $A$.\\
From \cite{Crisp-et-al} (or see \cite{Lothaire}) one deduces that $x=\Delta A$ is the fixed point of the morphism $\gamma$ given by
$$\gamma(1)=12, \quad \gamma(2)=121.$$
We define the extended morphism $\gamma_E$ on the alphabet $\{1,2,3,4\}$ by
$$\gamma_E(1)=13, \quad \gamma_E(2)=24, \quad \gamma_E(3)=241\quad \gamma_E(4)=132.$$
Note that $\gamma=\pi\gamma_E$, where $\pi(1)=\pi(2)=1$, and $\pi(3)=\pi(4)=2$. We define
$$y=1,3,2,4,1,2,4,1,3,2,1,3,\dots,$$
the fixed point of $\gamma_E$ with $y_1=1$. We claim that $y$ has the property that the letters 1 and 2 alternate in $y$. Indeed,  the words 132 and 12 are the only words in $y$ with prefix 1 and suffix 2 containing no 1's or 2's, and these are mapped  to
$$\gamma_E(12)=1324,\quad \gamma_E(132)=1324124,$$
in which 1's and 2's alternate,  and similarly the words 241 and 21 are mapped to 2413213 and 2413 in which 2's and 1's alternate. Since in the first case the first occurring letter is 1 and the last is 2, and in the second case the first occurring letter is 2 and the last is 1, it follows by induction that the letters 1 and 2 in $\gamma_E^n(1)$ alternate for all $n$.

We are interested in the positions 3,6,10,13,\dots\, of the letter 2 in $y$. Let $x'$ be defined by $x'_n=x_n-1$. Then $x'=0,1,0,1,0,0,1\dots$ is a Sturmian word with slope $\sqrt{2}-1$. Its mirror image $\tilde{x}'=1,0,1,0,1,1,0$ is a Sturmian word with slope $\tilde{\alpha}'=1-(\sqrt{2}-1)=2-\sqrt{2}$.  By Lemma 9.1.3 in \cite{all-shall}, the positions of 1's in $\tilde{x}'$ are given by the Beatty sequence $b=(\lfloor n \beta \rfloor)$, where
$$\beta =1/\tilde{\alpha}'= 1/(2-\sqrt{2}) = 1+\frac12\sqrt{2}.$$
But the 1's in $\tilde{x}'$ correspond to the 1's \emph{and} 2's in $y$, and since these alternate, the positions of the 2's in $y$ are given by the sequence $$(b_{2n})=(\lfloor 2n \beta \rfloor)=(\lfloor n(2+\sqrt{2})\rfloor).$$
Thus we found that the 2's in $y$ exactly occur at the Beatty complement $B$ of $A$.

\medskip

 { Step 2.} In this step we partition the `refinement' $y$ of the word $x=\Delta A$ in three words $w_1, w_2, w_3$, which will tell us how $\Delta AA$ behaves.  We claim that the three words
 $$w_1=132,\quad w_2=4,\quad w_3=124$$
 partition $y$. This follows directly from $\gamma_E(y)=y$ by noting that
 $$\gamma_E(w_1)=1324124=w_1w_2w_3,\quad   \gamma_E(w_2)=132=w_1,\quad  \gamma_E(w_3)=1324132=w_1w_2w_1.$$
 This equation induces a morphism $\sigma$ on the alphabet $\{1,\!2,3\}$, by replacing $w_j$ with $j$:
 $$\sigma(1)=123, \quad   \sigma(2)=1,  \quad  \sigma(3)=121.$$
How do we obtain $\Delta AA$ from $\Delta A$? Since $A(\N)=AA(\N)\cup AB(\N)$, a disjoint union,  one obtains $AA$ from $A$ by removing the integers $AB(n)$, which, of course, have index $B(n)$ in the sequence $A$. In Step 1 we showed that this sequence of indices corresponds to the positions of 2's in $y$. Now if such a 2 occurs in $w_1=132$, then the differences $x_k,x_{k+1},x_{k+2}=1,2,1$ in $x$ turn into differences 1,3 in $\Delta AA$, since  the second 1 disappears because of the removal of the $A$-number corresponding to $x_{k+2}$, and this 1 must be added to $x_{k+1}=2$.
The other possibility is that such a 2 occurs in $w_3=124$, and now the removal of the $A$-number corresponding to $x_{k+1}$ leads to differences 2,2 in  $\Delta AA$. The conclusion is that the decoration $\delta$ given by $\delta(1)=13, \delta(2)=2$ and $\delta(3)=22$ turns the fixed point of $\sigma$ into $\Delta AA$. \hfill $\Box$

\bigskip

\begin{corollary}\label{rem:morphic} Here is a way to write  $\Delta AA = 1,3,2,2,2,1,3,1,3,2,\dots$ as a morphic word (derived from the previous theorem). Let $\theta$ on $\{a,b,c,d\}$ be the morphism given by
$\theta: \; a\rightarrow adc,\, b\rightarrow adc, \,  c\rightarrow  ad, \, d\rightarrow bc.$\\
Let the letter-to-letter morphism $\lambda$ be given by
$\lambda: \;  a\rightarrow 1,\, b\rightarrow 2, \,  c\rightarrow  2, \, d\rightarrow 3.$
Then  $\Delta AA = \lambda( \theta^{\infty}(a)).$
\end{corollary}

\section*{Acknowledgement}
I am grateful to Jean-Paul Allouche for reading a draft of this paper, and suggesting that the $\alpha=\sqrt{2}$ case for the iterated Beatty sequences might be generalized to all metallic means.

\end{document}